\documentclass[a4paper, 12pt, leqno]{amsart}
\usepackage{mathrsfs}
\usepackage{amsmath}
\usepackage{amscd}
\usepackage{amssymb}
\usepackage{amsthm}
\usepackage{geometry}
\usepackage{tikz}

\title[Based Ring of Two-sided Cells ]
{The Based Rings of Two-sided cells in an Affine Weyl group of type $\tilde B_3$, I}

\author[Y. Qiu and N. Xi]{Yannan Qiu$^{*}$ and Nanhua XI$^{\dagger}$}
\address{$^{*}$
School of Mathematical Sciences\\
Zhejiang University \\
Zhejiang 310058, China \\
China} \email{qiuyannan@zju.edu.cn}
\address{$^{\dagger}$
Academy of Mathmatics and Systems Science\\
Chinese Academy of Sciences\\
Beijing 100190, China\\
and\\
School of Mathematical Sciences\\
University of Chinese Academy of Sciences\\
Chinese Academy of Sciences\\
Beijing 100049, China} \email{nanhua@math.ac.cn}
\thanks{N. Xi was partially supported by National Key R\&D Program of China, No. 2020YFA0712600, and by National Natural Science Foundation of China, No. 11688101.}

\begin{document}
\baselineskip=18pt
\begin{abstract}
For type $\tilde B_3$ we show that Lusztig's conjecture on the structure of the based ring of two-sided cell corresponding to the unipotent class
 in $Sp_6(\mathbb C)$ with 3 equal Jordan blocks needs modified.
\end{abstract}

\maketitle




We are concerned with the based ring of two-sided cells in an affine Weyl group of type $\tilde B_3$. In this paper we show Lusztig's conjecture on the structure of the based ring   of the two-sided cell  corresponding to the nilpotent
element in $Sp_6(\mathbb C)$ with 3 equal Jordan blocks needs modified, see section 4. This work was motivated by an inquiry of R. Bezrukanikov, who noted  (based on the work of Losev and Panin) that Lusztig's conjecture describing the summand of $J$ in terms of equivariant sheaves on the square of a finite set does not hold as stated for the 2-sided cell corresponding to the unipotent class in $Sp_6(\mathbb C)$ with 3 equal Jordan blocks and was interested  in explicit description of the summand in the asymptotic Hecke algebra $J$ attached to the affine Weyl group of type $\tilde B_3$. In a subsequent paper, we will discuss based ring of other two-sided cells in the affine Weyl group of type $\tilde B_3$.

The contents of the paper are as follows. Section 1 is devoted to preliminaries, which include some basic facts on (extended) affine Weyl groups and their Hecke algebras and formulation of  Lusztig's conjecture on the structure of based ring of a two-sided cell in an affine Weyl group. In section 2 we recall some results on cells of (extended) affine Weyl group of type $\tilde B_3$, which are due to J. Du. In section 3 we discuss the based ring $J_{\Gamma\cap\Gamma^{-1}}$ for a left cell in the concerned two-sided cell. In section 4 we show that Lusztig's conjecture on the structure of the based ring of the concerned two-sided cell needs modified.

After the paper is complete, we learned that R. Bezrukavnikov, S. Dawydiak and  G. Dobrovolska also worked out this example, see [BDD]. Their paper also contains discussion of relation of concerned central extensions to the cocenter of the affine Hecke algebra and the asymptotic affine Hecke algebra.

\def\Cal{\mathcal}
\def\bold{\mathbf}
\def\ca{\mathcal A}
\def\cdz{\mathcal D_0}
\def\cd{\mathcal D}
\def\cdo{\mathcal D_1}
\def\bold{\mathbf}
\def\l{\lambda}
\def\le{\leq}

\def\ll{\underset {L}{\leq}}
\def\rl{\underset {R}{\leq}}
\def\lr{\rl}
\def\lrl{\underset {LR}{\leq}}
\def\llr{\lrl}
\def\el{\underset {L}{\sim}}
\def\er{\underset {R}{\sim}}
\def\elr{\underset {LR}{\sim}}
\def\ds{\displaystyle\sum}

\section{Affine Weyl groups and their Hecke algebras}

In this section we recall some  basic facts on (extended) affine Weyl groups and their Hecke algebras and state Lusztig's conjecture on the structure of based ring of a two-sided cell  of an affine Weyl group.

\medskip

 \noindent{\bf 1.1. Extended affine Weyl groups and their Hecke algebras} \ \ Extended affine groups can be constructed through root systems or reductive groups. For our purpose, we use reductive groups to construct extended affine Weyl groups.

 Let $G$ be a connected reductive algebraic group over
the field {$\mathbb C$} of complex numbers  and $T$ a maximal torus of
$G$. Let $N_G(T)$ be the normalizer of $T$ in $G$. Then
$W_0=N_G(T)/T$ is a Weyl group, which acts on the character group
$X={\rm{Hom}}(T,\bold C^*)$ of $T$. The semi-direct product
$W_0\ltimes X$ is  called the extended affine Weyl group attached to $G$. The set of simple reflections of $W$ is denoted by $S$. We shall
denote the length function of $W$ by $l$ and use $\leq$ for the
Bruhat order on $W$. We also often write $y<w$ or $w>y$ if $y\le w$ and $y\ne w$.

Let $H$ be the  Hecke algebra  of $(W,S)$ over $\Cal A=\mathbb
Z[q^{\frac 12},q^{-\frac 12}]$ $(q$ an indeterminate) with
parameter $q$. Let  $\{T_w\}_{w\in W}$ be its standard basis. Then we have $(T_s-q)(T_s+1)=0$ and $T_wT_u=T_{wu}$ if $l(wu)=l(w)+l(u)$.
Let
$C_w=q^{-\frac {l(w)}2}\sum_{y\le w}P_{y,w}T_y,\ w\in W$ be the
Kazhdan-Lusztig basis of $H$, where $P_{y,w}$ are the
Kazhdan-Lusztig polynomials. The degree of $P_{y,w}$ is less than
or equal to $\frac12(l(w)-l(y)-1)$ if $y< w$ and  $P_{w,w}=1$. We write
$P_{y,w}=\mu(y,w)q^{\frac12(l(w)-l(y)-1)}$+lower degree terms   if $y< w$.

We shall write $y\prec w$ if $\mu(y,w)\ne 0$. The
coefficient $\mu(y,w)$ is very interesting, this can be seen even
from the recursive formulas below  for Kazhdan-Lusztig
polynomials ([KL]). Convention: set $P_{y,w}=0$ if $y\not\le w$.

\medskip
\def\vp{\varphi}
\def\st{\stackrel}
\def\sc{\scriptstyle}

(a) Let $y\le w$. Assume that $sw\le w$ for some $s\in S$. Then
$$P_{y,w}=q^{1-c}P_{sy,sw}+q^cP_{y,sw}-\sum_{\st {\st {z\in W}{y\le z\prec sw}}{sz<z}}\mu(z,sw)q^{\frac{l(w)-l(z)}2}P_{y,z},$$
where $c=1$ if $sy< y$ and $c=0$ if $sy>y$.

\medskip

(b) Let $y\le w$. Assume that $ws\le w$ for some $s\in S$. Then
$$P_{y,w}=q^{1-c}P_{ys,ws}+q^cP_{y,ws}-\sum_{\st {\st {z\in W}{y\le z\prec ws}}{zs<z}}q^{\frac{l(w)-l(z)}2}P_{y,z},$$
where $c=1$ if $sy< y$ and $c=0$ if $sy>y$.

\medskip

From the two formulas above one gets (see [KL])

\medskip

(c) Let $y,w\in W$ and $s\in S$ be such that $y<w,\ sw<w,\ sy>y$. Then $y\prec w$ if and only if $w=sy$. Moreover this implies that $\mu(y,w)=1$.

\medskip

(d) Let $y,w\in W$ and $s\in S$ be such that $y<w,\ ws<w,\ ys>y$. Then $y\prec w$ if and only if $w=ys$. Moreover this implies that $\mu(y,w)=1$.

\medskip

The following  formulas for computing $C_w$ (see [KL]) will be used in section 3 and 4.

\medskip

(e) Let $w\in W$ and $s\in S$, then
\begin{align} C_sC_w=\begin{cases}\displaystyle (q^{\frac12}+q^{-\frac12})C_w,\quad &\text{if\ }sw<w,\\
\displaystyle C_{sw}+\sum_{\st  {z\prec w}{sz<z}}\mu(z,w)C_z,\quad&\text{if\ }sw\ge w.\end{cases}\end{align}
\begin{align} C_wC_s=\begin{cases}\displaystyle (q^{\frac12}+q^{-\frac12})C_w,\quad &\text{if\ }ws<w,\\
\displaystyle C_{sw}+\sum_{\st  {z\prec w}{zs<z}}\mu(z,w)C_z,\quad&\text{if\ }ws\ge w.\end{cases}\end{align}

\medskip

\noindent{\bf 1.2. Cells of affine Weyl groups}  \  \ Cells in Coxeter groups were defined in [KL] through the coefficients $\mu(y,w)$. Cells in Coxeter groups can also be defined through Kazhdan-Lusztig basis. Following Lusztig [L1], we define preorders $\ll,\ \rl,\ \lrl$ on extended affine Weyl group $W$ as follows.

For $h,\, h'\in H$ and $x\in W$, write
\begin{alignat*}{2} hC_x&=\sum_{y\in W}a_yC_y,\qquad a_y\in \mathcal A,\\
 C_xh&=\sum_{y\in W}b_yC_y,\qquad b_y\in \mathcal A,\\
  hC_xh'&=\sum_{y\in W}c_yC_y,\qquad c_y\in \mathcal A.\end{alignat*}
  Define $y\ll x$ if $a_y\ne 0$ for some $h\in H$, $y\rl x$ if $b_y\ne 0$ for some $h\in H$, $y\lrl x$ if $c_y\ne 0$ for some $h,h'\in H$.

  We write $x\el y$ if $x\ll y\ll x$, and $x\er y$ if $x\rl y\rl x$, $x\elr y$ if $x\lrl y\lrl x$. Then $\el,\ \er,\ \elr$ are equivalence relations on $W$. The equivalence classes are called left cells, right cells, two-sided cells of $W$ respectively. Note that if $\Gamma$ is a left cell of $W$, then $\Gamma^{-1}=\{ w^{-1}\,|\, w\in\Gamma\}$ is a right cell.

\medskip

For $w\in W$, set $R(w)=\{s\in S\,|\, ws\le w\}$ and $L(w)=\{s\in S\,|\, sw\le w\}.$ Then we have (see [KL])

\medskip

(a) $R(w)\subset R(u)$ if $u\ll w$ and $L(w)\subset L(u)$ if $u\rl w.$ In particular, $R(w)= R(u)$ if $u\el w$ and $L(w)= L(u)$ if $u\er w.$

\medskip

\noindent{\bf 1.3.   $*$ operations}\ \ The $*$ operation introduced in [KL] and   generalized in [L1] is a useful tool in the theory of cells of Coxeter groups.

Let $s,t$ be simple reflections in $S$ and assume that $st$ has order  $m\ge 3$. Let $w\in W$ be such that $sw\ge w,\ tw\ge w$. The $m-1$ elements $sw,\ tsw,\ stsw,\ ...,$ is called a left string (with respect to $\{s,t\})$, and the $m-1$ elements $tw,\ stw,\ tstw,\ ..., $ is also called a left string (with respect to $\{s,t\})$. Similarly we define right strings (with respect to $\{s,t\}$). Then (see [L1])

\medskip

(a) A left string in $W$ is contained in a left cell of $W$ and a right string in $W$ is contained in a right cell of $W$.

\medskip
 Assume that $x$ is in a left (resp. right ) string (with respect to $\{s,t\})$ of length $m-1$ and is the $i$th element of the left (resp. right) string, we define ${}^*x$ (resp. $x^*$) to be the $(m-i)$th element of the string, where $*=\{s,t\}$. The following result is proved in [X2].

 \medskip

 (b) Let $x$ be in $W$ such that $x$ is in a left string with respect to $*=\{s,t\}$ and is also in a right string with respect to $\star=\{s',t'\}$. The ${}^*x$ is in a right string with respect to $\{s',t'\}$ and $x^\star$ is in a left string with respect to $\{s,t\}$. Moreover ${}^*(x^\star)=({}^*x)^\star$. We shall write ${}^*x^\star$
 for ${}^*(x^\star)=({}^*x)^\star$.

\medskip

The following result is due to Lusztig [L1].

\medskip

(c) Let $\Gamma$ be a left cell of $W$ and an element $x\in\Gamma$ is in a right string $\sigma_x$ with respect to $*=\{s,t\}$. Then any element $w\in\Gamma$ is in a right string $\sigma_w$ with respect to $*=\{s,t\}$. Moreover $\Gamma^*=\{w^*\,|\, w\in \Gamma\}$ is a left cell of $W$ and $\Omega=\displaystyle \left(\cup_{w\in\Gamma}\sigma_w\right)-\Gamma$ is a union of at most $m-2$ left cells.

\medskip

Following Lusztig [L1] we set $\tilde\mu(y,w)=\mu(y,w)$ if $y< w$ and $\tilde\mu(y,w)=\mu(w,y)$ if $w<y$. For convenience we also set $\tilde\mu(y,w)=0$ if $y\nless w$ and $w\nless y$. Assume that $x_1,x_2,...,x_{m-1}$  and $y_1,y_2,...,y_{m-1}$ are two left strings with respect to $*=\{s,t\}$. Define
$$a_{ij}=\begin{cases}\tilde\mu(x_i,y_j),\quad &\text{if\ } \{s,t\}\cap L(x_i)=\{s,t\}\cap L(y_j),\\
0,\quad &\text{otherwise}.\end{cases}$$
Lusztig proved the following identities (see subsection 10.4 in [L1]).

\medskip

(d) If $m=3$, then $a_{11}=a_{22}$ and $a_{12}=a_{21}$.

\medskip

(e) If $m=4$, then
\begin{align} a_{11}=a_{33},\ a_{13}=a_{31},\ a_{22}=a_{11}+a_{13},\ a_{12}=a_{21}=a_{23}=a_{32}.\end{align}

\medskip

\noindent{\bf 1.4. Lusztig's $a$-function}\quad
For $x,y\in W$, write $$C_xC_y=\sum_{z\in
W}h_{x,y,z}C_z,\qquad h_{x,y,z}\in \mathcal A
 =\mathbb Z[q^{\frac 12},q^{-\frac 12}].$$
 Following Lusztig ([L1]), we define
 $$a(z)={\rm min}\{i\in\bold N\ |\ q^{-\frac i2}h_{x,y,z}\in\mathbb Z[q^{-\frac
 12}]{\rm\ for \ all\ }x,y\in W\}.$$ If for any $i$,
 $q^{-\frac i2}h_{x,y,z}\not\in\mathbb Z[q^{-\frac
 12}]{\rm\ for \ some\ }x,y\in W$, we set $a(z)=\infty.$
 The following properties are proved in [L1].

 \medskip

 (a) Then $a(w)\le
l(w_0)$ for any $w\in W$, where $w_0$ is the longest element in the Weyl group $W_0$.

\medskip

(b) $a(x)\ge a(y)$ if $x\lrl y$. In particular, $a(x)=a(y)$ if $x\elr y$.

\medskip

(c) $x\el y$ (resp. $x\er y,\ x\elr y$) if $a(x)=a(y)$ and $x\ll y$ (resp. $x\lr y,\ x\llr y$).

Following Lusztig,
 we define  $\gamma_{x,y,z}$ by the following formula,
 $$h_{x,y,z}=\gamma_{x,y,z}q^{\frac {a(z)}2}+
 {\rm\ lower\ degree\ terms}.$$
Springer showed that $l(z)\ge a(z)$ (see [L2]). Let $\delta(z)$ be
the
 degree of $P_{e,z}$, where $e$ is the neutral element of $W$.
 Then actually one has $l(z)-a(z)-2\delta(z)\ge 0$ (see [L2]). Set
 $$\cd =\{z\in W\ |\ l(z)-a(z)-2\delta(z)=0\}.$$

 The elements of $\cd$ are involutions, called distinguished involutions of
 $(W,S)$ (see [L2]). The following properties are due to Lusztig [L2] except the trivial (i) and (j) [X2].

 (d) $\gamma_{x,y,z}\ne 0\Longrightarrow x\el y^{-1},\ y\el z,\ x\er z.$

 (e) $x\el y^{-1}$ if and only if $\gamma_{x,y,z}\ne 0$ for some $z\in W$.

 (f) $\gamma_{x,y,z}=\gamma_{y,z^{-1},x^{-1}}=\gamma_{z^{-1}, x,y^{-1}}$.

 (g) $\gamma_{x,d,x}=\gamma_{d,x^{-1},x^{-1}}=\gamma_{x^{-1},x,d}=1$ if $x\el d$ and $d$ is a distinguished involution.

 (h) $\gamma_{x,y,z}=\gamma_{y^{-1},x^{-1},z^{-1}}.$

 (i) If $\omega,\tau\in W$ has length 0, then
 $$\gamma_{\omega x,y\tau,\omega z\tau}=\gamma_{x,y,z},\ \
 \gamma_{ x\omega,\tau y, z}=\gamma_{x,\omega\tau y,z}.$$

 (j) Let $x,y,z\in W$ be such that (1) $x$ is in a left string with respect to $*=\{s,t\}$ and also in a right string with respect to $\#=\{s',t'\}$, (2) $y$ is in a left string with respect to $\#=\{s',t'\}$ and also  in a right string with respect to $\star=\{s'',t''\}$, (3) $z$ is in a left string with respect to $*=\{s,t\}$ and also  in a right string with respect to $\star=\{s'',t''\}$. Then
 $$\gamma_{x,y,z}=\gamma_{{}^*x^\#,{}^\#y^\star,{}^*z^\star}.$$

\medskip
\def\tt{\tilde T}

For $w\in W$, set $\tilde T_w=q^{-l(w)/2}T_w$. For $x,y\in W$, write
$$\tt_x\tt_y=\sum_{z\in W}f_{x,y,z}\tt_z,\qquad f_{x,y,z}\in\mathcal A=\mathbb Z[q^{\frac12},q^{-\frac12}].$$

(k) If $x,y,w$ are in a two-sided cell of $W$ and $f_{x,y,w}=\lambda q^{\frac{a(w)}2}+$ lower degree terms and as Laurent polynomials in $q^{\frac12}$, deg$f_{x,y,z}\le a(w)$ for all $z\in W$, then
$$\gamma_{x,y,w}=\lambda.$$

 (l) Each left cell (resp. each right cell) of $W$ contains a unique distinguished involution.

 (m) Each two-sided cell of $W$ contains only finitely many left cells.

\medskip

\def\ll{\underset {L}{\leq}}
\def\rl{\underset {R}{\leq}}

\def\lrl{\underset {LR}{\leq}}
\def\llr{\lrl}
\def\el{\underset {L}{\sim}}
\def\er{\underset {R}{\sim}}
\def\elr{\underset {LR}{\sim}}
\def\ds{\displaystyle\sum}

\def\vp{\varphi}
\def\st{\stackrel}
\def\sc{\scriptstyle}

 \noindent{\bf 1.5. The based ring of a two-sided cell}\quad For each two-sided cell $c$ of $W$, Let $J_c$ be the free $\mathbb Z$-module with a basis $t_w,\ w\in c$. Define
  $$t_xt_y=\sum_{z\in c}\gamma_{x,y,z}t_z.$$
  Then $J_c$ is an associative ring with unit $\sum_{d\in\cd\cap c}t_d.$

  The ring $J=\bigoplus_{c}J_c$ is a ring with unit $\sum_{d\in\cd}t_d$. Sometimes $J$ is called asymptotic Hecke algebra since Lusztig established an injective $\ca$-algebra homomorphism
  \begin{alignat*}{2}   \phi: H&\to J\otimes\ca,\\
  C_x&\to\sum_{\st {\st {d\in\cd}{w\in W}}{w\el d}}h_{x,d,w}t_w.\end{alignat*}

  \medskip

  \noindent{\bf 1.6. Lusztig's conjecture on the structure of $J_c$}\quad In [L3] Lusztig states a conjecture on $J_c$ using equivariant $K$-groups on finite sets.

  Assume that $F$ is a (possibly disconnected) reductive group over complex numbers field and $F$ acts algebraically on a finite set $Y$. Then an $F$-equivariant complex vector bundle (=$F$-v.b.) $V$ on $Y$ is just a collection of finite dimensional complex vector spaces $V_a,\ a\in Y$ with a given rational representation of $F$ on $\oplus_{a\in Y}V_a$ and $gV_a=V_{ga}$ for all $g\in F$ and $a\in Y$. The direct sum of $F$-v.b. can be defined naturally. Clearly the category of $F$-v.b. on $Y$ is abelian so that the Grothendieck group $K_F(Y)$ is well defined. The set $\mathcal T$ of  isomorphism classes of irreducible $F$-v.b. on $Y$ is a basis of $K_F(Y)$. For any $a\in Y$ let $F_a$ be the stabilizer of $a$ in $F$. It is easy to see that there exists a bijection between $\mathcal T$ and the set of pairs $(a,\rho)$ where $a\in Y$, $\rho\in\text{Irr}F_a$ (the set of isomorphism classes of irreducible rational representations of $F_a$), modulo the obvious action of $F$.

  Let $F$ act on $Y\times Y$ diagonally. Then $K_F(Y\times Y)$ is an associative ring if we define the multiplication $*$ in $K_F(Y\times Y)$ by setting
  $$(V*V')_{(a,b)}=\bigoplus_{r\in Y}V_{(a,r)}\otimes V'_{(r,b)}.$$
  The unit $\mathbf 1$ is the $F$-v.b. whose stalk at any point $(a,a)$ in the diagonal of $Y\times Y$ is the trivial representation of $F_a$ and is zero at any point off the diagonal of $Y\times Y$.

  Now let $G$ be a simply connected simple algebraic group over $\mathbb C$. Lusztig establishes a bijection between the two-sided cells of the extended affine Weyl group $W$ and the unipotent classes of $G$.

  For each two-sided cell $c$ of $W$, let $u$ be a unipotent element in the unipotent class corresponding to $c$ and let $F_c$ be a maximal reductive subgroup of the centralizer of $u$ in $G$.

  {\bf Conjecture} (Lusztig [L3]): There exists a finite set $Y$ with an algebraic action of $F_c$ and a  bijection
  $$\pi: c\to \text{isomorphism classes of irreducible $F_c$-v.b. on}\ Y\times Y.$$
  such that

  (i)  The bijection $\pi$ induces an ring homomorphism
  $$\pi: J_c\to K_{F_c}(Y\times Y),\ \ t_x\to \pi(x).$$

  (ii) $\pi(x^{-1})_{(a,b)}=\pi(x)_{(b,a)}^*$ is the dual representation of $\pi(x)_{(b,a)}.$

  \section{Cells in an extended affine Weyl group of type $\tilde B_3$}

  In this section $G=Sp_6(\mathbb C)$, so that the extended affine Weyl group $W$ attached to $G$ is of type $\tilde B_3$. The left cells and two-sided cells are described by J. Du (see [D]). We recall his results.

  {\bf 2.1. The Coxeter graph of $W$}. As usual, we number the 4 simple reflections $s_0,\ s_1,\ s_2,\ s_3$ in $W$ so that
  \begin{alignat*}{2} &s_0s_1=s_1s_0,\quad s_0s_3=s_3s_0,\quad s_1s_3=s_3s_1,\\
  &(s_0s_2)^3=(s_1s_2)^3=e,\quad (s_2s_3)^4=e,\end{alignat*}
  where $e$ is the neutral element in $W$. The relations among the simple reflections can be read through the following Coxeter graph:

 
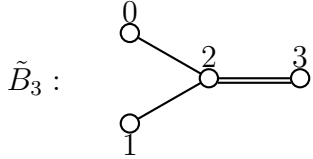
\begin{center}
  \begin{tikzpicture}[scale=.6]
    \draw (-1,0) node[anchor=east]  {$\tilde B_3:$};
    \draw[thick] (2 cm,0) circle (.2 cm) node [above] {$2$};
    \draw[xshift=2 cm,thick] (150:2) circle (.2 cm) node [above] {$0$};
    \draw[xshift=2 cm,thick] (210:2) circle (.2 cm) node [below] {$1$};
    \draw[thick] (4 cm,0) circle (.2 cm) node [above] {$3$};
    \draw[xshift=2 cm,thick] (150:0.2) -- (150:1.8);
    \draw[xshift=2 cm,thick] (210:0.2) -- (210:1.8);
    \draw[thick] (2.2,0) --+ (1.6,0);
    \draw[thick] (2.2,-0.1) --+ (1.6,0);
 \end{tikzpicture}
\end{center}

  There is a unique nontrivial element $\tau $ in $W$ with length 0. We have $\tau^2=e,\ \tau s_0\tau=s_1,\ \tau s_i\tau=s_i$ for $i=2,3.$ Note that $s_1,s_2,s_3$ generate the Weyl group $W_0$ and $s_0,s_1,s_2,s_3$ generate an affine Weyl group $W'$ of type $\tilde B_3$. And $W$ is generated by $\tau,\ s_0,s_1,s_2,s_3$.

  \medskip

  {\bf 2.2 Cells in $W$}\quad According to [D], the extended affine Weyl group $W$ attached to  $Sp_6(\mathbb C)$ has 8 two-sided cells:
    $$A,\quad B,\quad C,\quad D,\quad E, \quad F, \quad G,\quad H.$$
   The following table displays some useful information on these two-sided cells.
  \begin{center}
  \doublerulesep 0.4pt \tabcolsep8pt
\begin{tabular}{ccccc}
\hline
& &Number & Size of Jordan blocks& Maximal reductive subgroup\\
& & of left & of the  corresponding & of a unipotent element in\\
$X$ & a(X) & cells in $X$ &unipotent class in $Sp_6(\mathbb C)$&the corresp. unipotent class\\
\hline$A$&9&48&\hfill (111111)\ \ \ \ \ \ \ \ \ \ \ \   &$Sp_6(\mathbb C)$\\
$B$&6&24&\hfill (21111)\ \ \ \ \ \ \ \ \ \ \ \  &$Sp_4(\mathbb C)\times\mathbb Z/2\mathbb Z$\\
$C$&4&18&\hfill (2211)\ \ \ \ \ \ \ \ \ \ \ \  &$SL_2(\mathbb C)\times O_2(\mathbb C)$\\
$D$&3&12&\hfill (222)\ \ \ \ \ \ \ \ \ \ \ \  &$O_3(\mathbb C)$\\
$E$&2&8&\hfill (411)\ \ \ \ \ \ \ \ \ \ \ \  &$SL_2(\mathbb C)\times\mathbb Z/2\mathbb Z$\\
$F$&2&6&\hfill (33)\ \ \ \ \ \ \ \ \ \ \ \  &$SL_2(\mathbb C)$\\
$G$&1&4&\hfill (42)\ \ \ \ \ \ \ \ \ \ \ \  &$\mathbb Z/2\mathbb Z\times\mathbb Z/2\mathbb Z$\\
$H$&0&1&\hfill (6)\ \ \ \ \ \ \ \ \ \ \ \  &$\mathbb Z/2\mathbb Z$\\
\hline
\end{tabular}
\end{center}\vspace{2mm}

\medskip

{\bf 2.3} The $R$-set of a left cell $\Gamma$ is defined to be $R(\Gamma)=\{s\in S\,|\, xs\le x\}$ for any chosen $x\in \Gamma$. By 1.2(a), the $R$-set of a left cell is well defined. We use notation $X_{ij...}$ to indicate a left cell $X$ with $R$-set $\{s_i,s_j,...\}$.

 For a reduced expression  $s_{i_1}s_{i_2}\cdots s_{ik}$ of an element in $W$, we often write $i_1i_2\cdots i_k$ instead of the reduced expression.

 J. Du has given a representative for each left cell (see [D,Figure I, Theorem 6.4]. We recall his result for the two-sided cell $D$, which is the two-sided cell in $W$ corresponding to the unipotent class of $Sp_6(\mathbb C)$ with 3 equal Jordan blocks.

 (a) There are  12 left cells in the two-sided cell $D$ and a representative of each left cell in $ D$ are:
$$\begin{array}{lllllllll} &D_{013},&013;&D_2,&0132;&D_{02},&01320;& D_{12},&01321;\\
&&&&&&&&\\
 &D_{3},&01323;&D_{03},&013203; &D_{01},&013201; &D_{13},&013213; \\
&&&&&&&&\\
&D'_{2},&0132032; &\widehat{D'_{2}},&0132132; &D_{1},&01320321; &D_{0},&01321320.   \end{array}$$

\medskip

 The value of $a$-function on $D$ is 3.

 \medskip

 Let $\Gamma$ and $\Gamma'$ be two left cells of $W$. If $\Gamma'=\Gamma^*$ for some $*=\{s,t\}$ (see subsection 1.3 for definition of $*$ operation), then we write $\Gamma\ \overset{\{s, t\}}{\text{------}}\ \Gamma'$. The following result is easy to verify.

 \medskip

 {\bf Lemma 2.4.} Keep the notations in subsection 2.4. Then we have

 $$ D_{3}\ \overset{\{s_2, s_3\}}{\text{------}}\ D_{013}\ \overset{\{s_1, s_2\}}{\text{------}}\ D_{2}.$$

$$ \begin{aligned}
D_{0}\ \overset{\{s_0, s_2\}}{\text{------}}\ &\widehat{D'_{2}}\ \overset{\{s_2, s_3\}}{\text{------}}\ D_{12}\ \overset{\{s_0, s_2\}}{\text{------}}\ D_{01}\ \overset{\{s_1, s_2\}}{\text{------}}\ D_{02}\ \overset{\{s_2, s_3\}}{\text{------}}\ D'_{2}\ \overset{\{s_1, s_2\}}{\text{------}}\ D_{1}.\\
&\ | \{s_1, s_2\}\qquad\qquad\qquad\qquad\qquad\qquad\qquad\ |\{s_0, s_2\}\\
&D_{13}\qquad\qquad\qquad\qquad\qquad\qquad\qquad\qquad\  D_{03}
\end{aligned}$$

 \section{The based ring of the two-sided cell $D$, I}

 In this section we discuss the based ring of $D$, more precisely, we discuss the based ring $J_{\Gamma\cap\Gamma^{-1}}$ for a left cell $\Gamma$  in $D$. The unipotent class $\mathcal C_D$ in $G=Sp_6(\mathbb C)$ corresponding to $D$ has three equal Jordan blocks, so a maximal reductive subgroup of $C_G(u)$ for any $u\in \mathcal C_D$ is isomorphic to $O_3(\mathbb C)$.  From now on, all representations are rational representations.

 \def\otc{O_3(\mathbb C)}
 \def\sotc{SO_3(\mathbb C)}

 Note thta $\otc=\mathbb Z/2\mathbb Z\times\sotc$. Denote by $\epsilon$ the sign representation of $\mathbb Z/2\mathbb Z$ and by $V(2k)$ the irreducible representation of $\sotc$ with highest weight $2k\ (k\ge 0)$. Both $\epsilon$ and $V(2k)$ are naturally irreducible representations of $\otc$. Up to isomorphism, all irreducible representations are $V(2k),\ \epsilon\otimes V(2k).$ We shall write $\epsilon V(2k)$ for $\epsilon\otimes V(2k),$ $k=0,1,2,3,....$

 \medskip
 \def\otc{O_3(\mathbb C)}
 \def\sotc{SO_3(\mathbb C)}

 {\bf Theorem 3.1} Let $\Gamma$ be a left cell in $D$. There is a natural bijection
 $$\pi: \Gamma\cap \Gamma^{-1}\to \text{Irr}O_3(\mathbb C)$$
 and the bijection induces a ring isomorphism
 $$J_{\Gamma\cap\Gamma^{-1}}\to \text{Rep}O_3(\mathbb C),\quad t_x\to \pi(x),$$
 where Irr$O_3(\mathbb C)$ is the set of isomorphism classes of irreducible   representations of $O_3(\mathbb C)$, Rep$\otc$ is the Grothendieck group of the category of   representations of $\otc$, the multiplication in Rep$\otc$ is given by tensor product of representations.

 {\bf Proof.} It is easy to see that if $\Gamma'=\Gamma^*$ for some $*=\{s,t\}$ with $s,t\in S$, then the map $\Gamma\cap\Gamma^{-1}\to\Gamma'\cap\Gamma'^{-1},\ x\to{}^*x^*$ is a bijection. Using   1.4(j), we see that the map $J_{\Gamma\cap\Gamma^{-1}}\to J_{\Gamma'\cap\Gamma'^{-1}},\ t_x\to t_{{}^*x^*}$ is a ring isomorphism.

 By Lemma 2.4, it suffices to prove Theorem 3.1  for $\Gamma=D_{013}$ and $\Gamma=D_{12}$. By Lemma 2.4, $D_{12}$ can be obtained from $D_0$ by two $*$-operations. According to [B], for $D_0$, Theorem 3.1  is true. By 1.4(j), Theorem 3.1  is true for $D_{12}$. So we only need to prove Theorem 3.1  for $\Gamma=D_{013}$. Also, we offer here a different proof for $D_{12}$ which exhibits the bijection $\pi:D_{12}\cap D_{12}^{-1}\to \text{Irr}O_3(\mathbb C)$ explicitly.

 The following two subsections are devoted to prove Theorem 3.1  for $D_{12}$ and $D_{013}$ respectively.

 \medskip

 {\bf 3.2.} \quad In this subsection we prove   Theorem 3.1 for $\Gamma=D_{12}$ by explicitly constructing the map $\pi$.  According to [D, Theorem 6.4] we have
$$ D_{12}\cap D_{12}^{-1}=\{ (s_1s_2s_3s_0)^ks_1s_2s_1, \  \tau s_2s_0s_1s_2s_1,\  \tau s_0s_2s_3s_0(s_1s_2s_3s_0)^ks_1s_2s_1\,|\, k\ge 0\}.$$

{\bf Claim 1:} The bijection
\begin{alignat*}{2} \pi:  D_{12}\cap D_{12}^{-1}&\longrightarrow \text{Irr}O_3(\mathbb C),\\
(s_1s_2s_3s_0)^ks_1s_2s_1&\longrightarrow V(2k),\\
 \tau s_2s_0s_1s_2s_1&\longrightarrow\epsilon,\\
   \tau s_0s_2s_3s_0(s_1s_2s_3s_0)^ks_1s_2s_1&\longrightarrow\epsilon V(2k+2),\end{alignat*}
   induces a ring isomorphism
    $$J_{D_{12}\cap D_{12}^{-1}}\longrightarrow \text{Rep}O_3(\mathbb C),\quad t_x\longrightarrow \pi(x).$$

    Now we prove the claim. Note that for any $x$ in $D_{12}\cap D_{12}^{-1}$, we have $x=x^{-1}$. Using 1.4 (h) we get

    \medskip

    (a) $t_xt_y=t_yt_x$ for any $x,y\in D_{12}\cap D_{12}^{-1}$.

    \medskip

    Let
    $$x_k=(s_1s_2s_3s_0)^ks_1s_2s_1,\ x'_0=\tau s_2s_0s_1s_2s_1, \ x'_{k+1}=\tau s_0s_2s_3s_0(s_1s_2s_3s_0)^ks_1s_2s_1.$$
    Thanks to identity in (a) above, to see Claim 1 we only need to verify the following identities:
    \begin{align} t_{x_k}t_{x_l}&=\sum_{|k-l|\le m\le|k+l|}t_{x_m},\\
    t_{x'_k}t_{x'_l}&=\sum_{|k-l|\le m\le|k+l|}t_{x_m},\\
    t_{x_k}t_{x'_l}&=\sum_{|k-l|\le m\le|k+l|}t_{x'_m}.\end{align}

     The identities (4), (5)  (6) are equivalent. In fact, we have the following result which  is easy to verify.

    \medskip

    (b) Let $*=\{s_2,s_0\}$ and $\#=\{s_1,s_2\}$. Then for any nonnegative integer $k$ we have
    $$x'_k=\tau({}^\#({}^*x_k))=((x_k^*)^\#)\tau.$$

     Note that $\tau^2=e$ is the neutral element in $W$ and $\tau s_0=s_1\tau,\ \tau s_2=s_2\tau$. Thus $(z\tau)^\#=(z^*)\tau$ and $(z\tau)^*=(z^\#)\tau$ for any $z\in D_{12}\cap D_{12}^{-1}$.  Using (b), 1.4(i) and 1.4(j), we get
     \begin{alignat*}{2} \gamma_{x'_k,x'_l,z}&=\gamma_{((x_k^*)^\#)\tau,\tau({}^\#({}^*x_l)),z}\\
     &=\gamma_{(x_k^*)^\#,\tau^2({}^\#({}^*x_l)),z}=\gamma_{(x_k^*)^\#,{}^\#({}^*x_l),z}\\
     &=\gamma_{x_k^*,{}^*x_l,z}=\gamma_{x_k,x_l,z};\\
      \gamma_{x_k,x'_l,z}&=\gamma_{x_k,((x_l^*)^\#)\tau,z}\\
     &=\gamma_{x_k,(x_l^*)^\#,z\tau}=\gamma_{x_k,x_l^*,(z\tau)^\#}\\
     &=\gamma_{x_k,x_l^*,(z^*)\tau}=\gamma_{x_k,x_l,((z^*)\tau)^*}\\
&=\gamma_{x_k,x_l ,((z^*)^\#)\tau}.\end{alignat*}
Therefore, the identities (4), (5) and (6) are equivalent. Now we prove identity (4). When $k=0$, $x_k=s_1s_2s_1$ is a distinguished involution and the identity (4) is trivial in this case. Similarly, identity (4) is true if $l=0$.

Now assume that $k=1$ and $l\ge 1$.
Let $\zeta=q^{\frac12}-q^{-\frac12}$. By a simple computation we get
$$\tilde T_{x_1}\tilde T_{x_1}=\zeta^3(\tilde T_{x_{2}}+\tilde T_{x_1}+\tilde T_{x_{0}})+\text{lower degree terms}.$$
Using 1.4(k) we know that identity (4) is true for $k=l=1.$

For $l\geq 2$, again by direct computation we get
\begin{alignat*}{2} \tilde T_{x_1}\tilde T_{x_l}&=\zeta^3(\tilde T_{x_{l+1}}+\tilde T_{x_l}+\tilde T_{x_{l-1}}+\tilde T_{w_{012}s_3s_2s_1(s_0s_3s_2s_1)^{l-2}}+\tilde T_{s_2s_1w_{023}s_1s_2(s_0s_3s_2s_1)^{l-2}})\\
&\qquad +\text{lower degree terms},\end{alignat*}
where $w_{012}$ is the longest element in the parabolic subgroup of $W$ generated by $s_0,s_1,s_2$ and $w_{023}$ is the longest element in the parabolic subgroup of $W$ generated by $s_0,s_2,s_3$. Clearly
\begin{alignat*}{2}&a(w_{012}s_3s_2s_1(s_0s_3s_2s_1)^{l-2})\geq a(w_{012})=6, \\ &a(s_2s_1w_{023}s_1s_2(s_0s_3s_2s_1)^{l-2})=9.\end{alignat*}
 So $w_{012}s_3s_2s_1(s_0s_3s_2s_1)^{l-2}$ and $s_2s_1w_{023}s_1s_2(s_0s_3s_2s_1)^{l-2}$ are not in $D$.

Using 1.4(k) we see that the identity (4) is true if $k=1,\ l\ge 2.$ Thus identity (4) is true for $k=1$ and all nonnegative integer $l$.

Assume that $k\ge 2$. What we have just proved tells us that $t_{x_k}=t_{x_1}t_{x_{k-1}}-t_{x_{k-1}}-t_{x_{k-2}}$. This means that we can use induction on $k$ to prove identity (4). The argument for identity (4) is complete and Claim 1 follows.

\medskip

 {\bf 3.3.} \quad In this subsection we prove  Theorem 3.1 for $\Gamma=D_{013}$ by explicitly constructing the map $\pi$.  According to [D, Theorem 6.4] we have
$$ D_{013}\cap D_{013}^{-1}=\{ (s_0s_1s_3s_2)^ks_0s_1s_3, \quad  \tau (s_0s_1s_3s_2)^ks_0s_1s_3\,|\, k\ge 0\}.$$

{\bf Claim 2:} The bijection
\begin{alignat*}{2} \pi:  D_{013}\cap D_{013}^{-1}&\longrightarrow \text{Irr}O_3(\mathbb C),\\
(s_0s_1s_3s_2)^ks_0s_1s_3&\longrightarrow V(2k),\\
 \tau (s_0s_1s_3s_2)^ks_0s_1s_3&\longrightarrow\epsilon V(2k),\end{alignat*}
   induces a ring isomorphism
    $$J_{D_{013}\cap D_{013}^{-1}}\longrightarrow \text{Rep}O_3(\mathbb C),\quad t_u\longrightarrow \pi(u).$$

    Now we prove claim 2. Note that for any $u$ in $D_{013}\cap D_{013}^{-1}$, we have $u=u^{-1}$. Using 1.4 (h) we get

    \medskip

    (c) $t_ut_v=t_vt_u$ for any $u,v\in D_{013}\cap D_{013}^{-1}$.

    \medskip

    Let
    $$u_k=(s_0s_1s_3s_2)^ks_0s_1s_3,\ \ u'_k=\tau (s_0s_1s_3s_2)^ks_0s_1s_3.$$
    Thanks to identity in (c) above, to see Claim 2 we only need to verify the following identities:
    \begin{align} t_{u_k}t_{u_l}&=\sum_{|k-l|\le m\le|k+l|}t_{u_m},\\
    t_{u'_k}t_{u'_l}&=\sum_{|k-l|\le m\le|k+l|}t_{u_m},\\
    t_{u_k}t_{u'_l}&=\sum_{|k-l|\le m\le|k+l|}t_{x'_m}.\end{align}

      Using 1.4 (i) we see that the  identities (7), (8),  (9) are equivalent.
      Now we prove identity (7). When $k=0$, $u_k=s_0s_1s_3$ is a distinguished involution and the identity (4) is trivial in this case. Similarly, identity (7) is true if $l=0$.

Now assume that $k=1$ and $l\ge 1$. Let $\xi=q^{\frac12}+q^{-\frac12}$. By a simple computation we get
\begin{align} C_{s_0s_1s_3s_2 s_0s_1s_3}&=C_{s_0}C_{s_1}C_{s_3}C_{s_2}C_{s_0s_1s_3}-(\xi^2+1)C_{s_0s_1s_3},\\
C_{s_0s_1s_3}C_{u_l}&=\xi^3C_{u_l}.\end{align}
Hence
\begin{align}  C_{u_1}C_{u_l}=C_{s_0s_1s_3s_2 s_0s_1s_3}C_{u_l}=\xi^3C_{s_0}C_{s_1}C_{s_3}C_{s_2}C_{u_l}-\xi^3(\xi^2+1)C_{u_l}.\end{align}

We compute $C_{s_0}C_{s_1}C_{s_3}C_{s_2}C_{u_l}$ step by step. Before continuing, we make a convention: {\sl we shall use the symbol $\Box$ for any element in the $\mathcal A$-submodule $H^{<013}$ of $H$ spanned by all $C_w$ with $a(w)\ge 4$. Then $\Box+\Box=\Box$ and $h\Box=\Box$ for any $h\in H$.}

\medskip

{\bf Step 1.} Compute $C_{s_2}C_{u_l}.$  By formula (1) in 1.1 (e), we have
$$C_{s_2}C_{u_l}=C_{s_2u_l}+\sum\limits_{\substack{z\prec u_l\\ s_2z<z}}\mu(z, u_l)C_z.$$

Note that $L(u_l)=\{s_0,s_1,s_3\}$. Assume that $z\prec u_l$ and $s_2z\le z$. If $s_3z\le z$, then $\{s_2,s_3\}\subset L(z)$ and $a(z)\ge a(s_2s_3s_2s_3)=4$. In this case, we have $C_z\in H^{<013}$. If $s_3z\ge z$,   we must have $s_3z=u_l$. Then $s_2z=s_2s_3u_l\ge s_3u_l=z$. This contradicts $s_2z\le z$. Therefore we  have
\begin{align} C_{s_2}C_{u_l}=C_{s_2u_l}+\Box\in C_{s_2u_l}+H^{<013}.\end{align}

\medskip

{\bf Step 2.}  Compute $C_{s_3}*C_{s_2u_l}$.

We have
$$C_{s_3}*C_{s_2u_l}=C_{s_3s_2u_l}+\sum\limits_{\substack{z\prec s_2u_l\\ s_3z<z}}\mu(z, s_2u_l)C_z.$$

Note that $ {L}(s_2u_l)=\{s_2\}$. Assume that $z\prec s_2u_l$ and $s_3z\le z$. If $s_2z\le z$, then $\{s_2,s_3\}\subset L(z)$ and $a(z)\ge a(s_2s_3s_2s_3)=4$. In this case, we have $C_z\in H^{<013}$. If $s_2z\ge z$,   we must have $s_2z=s_2u_l$. Then $ z= u_l$. Hence we have
\begin{align} C_{s_3}C_{s_2u_l}=C_{s_3s_2u_l}+C_{u_l}+\Box\in C_{s_3s_2u_l}+C_{u_l}+H^{<013}.\end{align}

\medskip

{\bf Step 3.} Compute $C_{s_1}*C_{s_3s_2u_l}$.

We have
$$C_{s_1}*C_{s_3s_2u_l}=C_{s_1s_3s_2u_l}+\sum\limits_{\substack{z\prec s_3s_2u_l\\ s_1z<z}}\mu(z, s_3s_2u_l)C_z.$$

Note that $ {L}(s_3s_2u_l)=\{s_3\}$. Assume that $z\prec s_3s_2u_l$ and $s_1z\le z$.

If $s_3z\ge z$, then $s_3z=s_3s_2u_1$. Hence $z=s_2u_1$ and $s_1z\ge z$, which contradicts the assumption $s_1z\le z$.

Now assume  $s_3z\le z$. Then $z=s_3s_1z'$ for some $z'\in W$ with $s_1z'\ge z',\ s_3z'\ge z'$.
When $s_2s_1z'\le s_1z'$, we have $z=s_3s_2s_1s_2\tilde z$ with $l(z)=4+l(\tilde z)$. If further $s_3s_1s_2\tilde z\le s_1s_2\tilde z$, then $s_1s_2\tilde z=s_1s_2s_3s_2s_3 y$ with $l(s_1s_2\tilde z)=5+l(y)$. In this case $a(z)\ge a(s_2s_3s_2s_3)=4$ and $C_z\in H^{<013}$.   If $s_3s_1s_2\tilde z\ge s_1s_2\tilde z$, using formula (3) in 1.3(e) we get
$$\mu(z,s_3s_2u_l)=\tilde \mu(z,u_l).$$
\begin{itemize}
\item[*] If $s_0z<z<u_l$, then $s_0\in {L}(s_2s_1s_2\tilde z)$. Hence   $a(z)\geq a(s_0s_1s_2s_0s_1s_2)=6$ and  $C_z\in H^{<013}$.
\item[*] If $s_0z>z$ and $z<u_l$, then $s_0z=u_l$. Hence  $z=s_1s_3s_2u_{l-1}$.
\item[*] If $\tilde\mu(y, v_k)=\mu(v_k, y)$, then we have $u_l\prec y\prec  s_3s_2u_l$, which forces $z=s_2u_l$. It contradicts  $s_3z\le z$.
\end{itemize}

When $s_2s_1z'\ge s_1z'$, using formula (3) in 1.3(e) we get
$$\mu(z,s_3s_2u_l)=\tilde\mu (s_3s_2z,u_l).$$
\begin{itemize}
\item[*] If $s_3s_2z\prec u_l$, and $s_1s_3s_2z\le z$, using  $z= s_1s_3z'$ we see that $s_2s_3z'\le s_3z'$. This forces that  $a(z)\geq a(s_2s_3s_2s_3)=4$. In this case we have   $C_z\in H^{<013}$.
\item[*] If $s_3s_2z\prec u_l$, and $s_1s_3s_2z\ge z$, then  we have  $s_1s_3s_2z=u_l$. This implies that $z=s_0s_2(s_1s_3s_2s_0)^{l-1}s_1s_3$.  It contradicts assumption $s_3z\le z$.
\item[*] If $u_l\prec s_3s_2z$, and $z\prec s_3s_2 u_l$, then $l(z)+2=l(u_l)+1$ or $l(z)+1=l(u_l)+2$. Hence $u_l=s_2z$ or $z=s_2u_l$. But then $s_1z\ge z$, which contradicts assumption $s_1z\le z$.
\end{itemize}

In conclusion, if $z\prec s_3s_2u_l$ and $s_1z\le z$, then either $C_z\in H^{<013}$ or $z=s_1s_3s_2u_{l-1}$. Hence we have
\begin{align} C_{s_1}C_{s_3s_2u_l}= C_{s_1s_3s_2u_l}+C_{s_1s_3s_2u_{l-1}}+\Box\in C_{s_1s_3s_2u_l}+C_{s_1s_3s_2u_{l-1}}+H^{<013}.\end{align}

\medskip

{\bf Step 4.} Compute $C_{s_0}C_{s_1s_3s_2u_l}$ and $C_{s_0}C_{s_1s_3s_2u_{l-1}}$. We only need to compute the first one. We have
$C_{s_0}C_{s_1s_3s_2u_l}=C_{u_{l+1}}+\sum\limits_{\substack{z\prec s_1s_3s_2u_l\\ s_0z<z}}\mu(z, s_1s_3s_2u_l)C_z$.

Assume that $ z\prec s_1s_3s_2u_l$ and $s_0z\le z$. Note that $ {L}(s_1s_3s_2u_l)=\{s_1,s_3\}$.

 If $s_1z\ge z$, we have  $s_1z=s_1s_3s_2u_l$. Then $z=s_3s_2u_l$ and  $s_0z\ge z$, which contradicts assumption $s_0z\le z$.

 Similarly, if $s_3z\ge z$, we have  $s_3z=s_1s_3s_2u_l$ and $z=s_1s_2u_l$. Then $s_0z\ge z$, which contradicts assumption $s_0z\le z$.

 Now assume  $s_1z\le z,\ s_3z\le z$. Then there exists $z'\in W$ such that   $z=s_0s_1s_3z'=s_3s_0s_1z'$ with $s_iz'\ge z$ for $i=0,1,3.$

 If $s_2s_0s_1z'\le s_0s_1z'$, then $a(z)\geq a(s_0s_1z')\geq a(s_0s_1s_2s_0s_1s_2)=6$. Hence $C_z\in H^{<013}$.

 If $s_2s_0s_1z'\ge s_0s_1z'$, using formula (3) in 1.3(e) we get $\mu(z, s_1s_3s_2u_l)=\tilde\mu(s_2z, s_1s_2u_l)$.
\begin{itemize}
\item[*] When $s_2z\prec s_1s_2u_l$, and $s_1s_2z\le s_2z$, we have $s_2s_3s_0z'<s_3s_0z'$. This implies that  $a(z)\geq a(s_2s_3s_2s_3)=4$, so $C_z\in H^{<013}$ in this case.
\item[*]  When $s_2z\prec s_1s_2u_l$, and $s_1s_2z\ge s_2z$, we have  $s_1s_2z=s_1s_2u_l$. Hence $z=u_l$.
\item[*] If $s_1s_2u_l\prec s_2z$, and $z\prec s_1s_3s_2u_{l}$, then $l(u_l)+3=l(z)+1$ and $s_2z=s_1s_3s_2u_{l}$. But $s_2s_1s_3s_2u_{l}\ge s_1s_3s_2u_{l}$, so $s_2z=s_1s_3s_2u_{l}$ could not occur.
\end{itemize}

In conclusion, if $ z\prec s_1s_3s_2u_l$ and $s_0z\le z$, then either $C_z\in H^{<013}$ or $z=u_{l }$. Hence we have
\begin{align} C_{s_0}C_{s_1s_3s_2u_l}&= C_{u_{l+1}}+C_{ u_{l}}+\Box,\\
C_{s_0}C_{s_1s_3s_2u_{l-1}}&= C_{u_{l }}+C_{ u_{l-1}}+\Box.\end{align}
(The second one follows from the first one.)

Note that $C_{s_0}C_{s_1}C_{u_l}=\xi^2 C_{u_l}$. Combining formulas (12)-(17) we get
$$\begin{aligned}
C_{u_1}C_{u_l}&= \xi^3(C_{u_{l+1}}+C_{u_l}+C_{u_{l}}+C_{u_{l-1}}+\xi^2C_{u_l})-\xi^3(\xi^2+1)C_{u_l}+\Box\\
&=\xi^3(C_{u_{l+1}}+C_{u_l}+ C_{u_{l-1}})+\Box \in \xi^3(C_{u_{l+1}}+C_{u_l}+ C_{u_{l-1}})+ H^{<013},
\end{aligned}$$

 The argument for identity (7) is complete and Claim 2 follows.

 We have completed the proof for Theorem 3.1.

\section{The based ring of the two-sided cell $D$, II}

 In this section we show that Lusztig's conjecture on the structure of $J_D$ needs a modification. For simplicity  we shall write $F$ for $\otc$.

 \medskip

 {\bf 4.1.} Assume there existed a finite $F$-set $Y$ and a bijection

 \centerline{$\pi:D\to$ the set of isomorphism classes of irreducible $F$-v.b. on $Y\times Y$}

 \noindent such that $\pi(x^{-1})_{(a,b)}=\pi(x)_{(b,a)}^*$ is the dual representation of $\pi(x)_{(b,a)}$ and the map $t_w\to \pi(w)$ defines a ring isomorphism (preserving the unit element) between $J_D$ and $K_F(Y\times Y)$.

 \medskip

 {\bf 4.2.} Keep the assumption in 4.1. Let $w\in D$, then $\pi(w)$ is an isomorphism class of irreducible $F$-v.b on $Y\times Y$. Note that the  irreducible $F$-v.b. on $Y\times Y$ is one to one corresponding to the set of triples $(x,y,\rho)$, where $(x,y)\in Y\times Y$ and $\rho$ is an irreducible   representation of the stabilizer $F_{x,y}$ of $(x,y)$, modulo the obvious action of $F$.

 Since $F=\mathbb Z/2\mathbb Z\times SO_3(\mathbb C)$, and $\sotc$ is connected so it acts on $Y\times Y$ trivially, each $F_{x,y}$ is either $\otc$ or $\sotc$. Then  an irreducible representation $\rho$ of $F_{x,y}$ essentially is  an irreducible representation of $\sotc$, at most differing from a sign representation of $\mathbb Z/2\mathbb Z$. Thus $\rho$ has a highest weight, which is a non-negative even integer $2k$ and the dimension of $\rho$ is $2k+1$. We will need the following two simple facts.

 \medskip

 (a) Let  $\rho$ and $\phi$ be two irreducible   representations of $\otc$ or of $\sotc$. Then the tensor product $\rho\otimes\phi$ can be decomposed into a direct sum of   $\min\{\dim\rho,\ \dim\phi\}$ irreducible representations of $\otc$ or of $\sotc$. In particular, the number of composition factors of $\rho\otimes\phi$ is odd.

 \medskip

 (b) Any   representation $\rho$  of $\otc$ or of $\sotc$ is isomorphic to its dual representation $\rho^*$.

 \medskip

 Let $w=s_0s_1s_3s_2s_1\in D_{12}$ (for notations see subsection 2.3). Then $w^{-1}=s_1s_2s_3s_1s_0\in D_{013}$. We first compute $t_{w^{-1}}t_w$ in $J$. To do this, we need to compute $C_{w^{-1}}C_w$ in the affine Hecke algebra. For simplicity, {\bf we often write $C_{i_0i_2\cdots i_m}$ for $C_x$ if $s_{i_0}s_{i_2}\cdots s_{i_m}$ is a reduced expression of $x$.}
 Note that $s_1,s_0,s_3$ are mutually commutative and $(s_0s_2)^3=(s_1s_2)^3=(s_2s_3)^4=e$ is the neutral element in $W$. In the affine Hecke algebra $H$ of $W$ we also use $e$ for the unit $C_e$.

We have
 \begin{alignat*}{2} C_{w^{-1}}&=C_{12310}=C_1C_2C_{310}-C_{310},\\
 C_{w^{-1}}C_w&=(C_1C_2-e)C_{310}C_{01321}\\
 &=\xi^3(C_1C_2C_{01321}-C_{01321})\\
 &=\xi^3(C_1C_{201321}-C_{01321})\\
 &\ \ \text{(note that \ }\mu(121,201321)=\mu(0121,01321)=1.)\\
 &=\xi^3(C_{1210321}+C_{01321}+C_{121}-C_{01321})\\
 &=\xi^3(C_{1210321}+C_{121}),\end{alignat*}
 where $\xi=q^{\frac12}+q^{-\frac12}$. We then get
 \begin{align} t_{w^{-1}}t_w=t_{1210321}+t_{121},\end{align}
 here we use notation $t_{i_1i_2\cdots i_m}$ for $t_x$ if $s_{i_1}s_{i_2}\cdots s_{i_m}$ is a reduced expression of $x$.

Note that the following formula is just formula (4) in subsection 3.2 for $k=l=1$.
\begin{align} t_{1210321}t_{1210321}=t_{12103210321}+t_{1210321}+t_{121}.\end{align}

\medskip

{\bf 4.3.} Keep the assumption in 4.1. For $w^{-1}=s_1s_2s_3s_1s_0$, let $\pi(w^{-1})=V$ be the irreducible $F$-v.b. on $Y\times Y$ corresponding to the class containing $(x,y,\rho)$. Then $\pi(w )=V'$ is the irreducible $F$-v.b. on $Y\times Y$ corresponding to the class containing $(y,x,\rho^*)$. So  $V_{x,y}=\rho$ is an irreducible   representation of $F_{x,y}$ and $V'_{y,x}=\rho^*=\rho$.

There are four possibilities for the  $F$-orbit containing $(x,y)$.

\medskip

Case 1.  Both $x$ and $y$ are $F$-stable. In this case $\rho$ is an irreducible representation of $\otc$ and $(V*V')_{p,q}\ne 0$ if and only if $(p,q)=(x,x)$. Moreover
$$(V*V')_{x,x}=V_{x,y}\otimes V'_{y,x}=\rho\otimes\rho.$$
According to fact 4.2(a), $\rho\otimes\rho$ can be decomposed into direct sum of  $\dim\rho$ irreducible representations of $\otc$. This means $V*V'$ can be decomposed into a direct sum of  $\dim\rho$ irreducible $F$-vector bundles on $Y\times Y$. By assumption 4.1 and fact 4.2(a), for $\displaystyle t_{w^{-1}}t_w=\sum_z\gamma_{w^{-1},w,u}t_u$, we get that  $\sum_u\gamma_{w^{-1},w,u}=\dim\rho$ is  odd. This contradicts to formula (18) in subsection 4.2, which shows $\sum_u\gamma_{w^{-1},w,u}=2.$ So this case would not occur.

\medskip

Case 2.  Both $x$ and $y$ are not $F$-stable. Then both $Fx$ and $Fy$ have cardinality 2. Let $Fx=\{x,p\},\ Fy=\{y,q\}$. In this case $F_{x,y}=SO_3(\mathbb C)$ and $\rho$ is an irreducible representation of $\sotc$. We have
 $V_{\alpha,\beta}=0$ if $(\alpha,\beta)\ne (x,y),\ (p,q)$, and $V'_{\alpha,\beta}=0$ if $(\alpha,\beta)\ne (y,x),\ (q,p)$. Thus  $(V*V')_{\alpha,\beta}=0$ if $(\alpha,\beta)\ne (x,x),\ (p,p)$. Moreover
$$(V*V')_{x,x}=V_{x,y}\otimes V'_{y,x}=\rho\otimes\rho*=\rho\otimes\rho.$$
According to fact 4.2(a), $\rho\otimes\rho$ can be decomposed into direct sum of  $\dim\rho$ irreducible representations of $\sotc$. This means $V*V'$ can be decomposed into direct sum of  $\dim\rho$ irreducible $F$-vector bundles on $Y\times Y$. By assumption 4.1 and fact 4.2(a), for $\displaystyle t_{w^{-1}}t_w=\sum_z\gamma_{w^{-1},w,u}t_u$, we get that  $\sum_u\gamma_{w^{-1},w,u}=\dim\rho$ is  odd. This contradicts to formula (18) in subsection 4.2, which shows $\sum_u\gamma_{w^{-1},w,u}=2.$ So this case would not occur.

\medskip

Case 3. The element $x$ is not $F$-stable and $y$ is $F$-stable. Then  $Fx$ has cardinality 2. Let $Fx=\{x,p\}$. In this case $F_{x,y}=SO_3(\mathbb C)$ and $\rho$ is an irreducible representation of $\sotc$. We have
 $V_{\alpha,\beta}=0$ if $(\alpha,\beta)\ne (x,y),\ (p,y)$, and $V'_{\alpha,\beta}=0$ if $(\alpha,\beta)\ne (y,x),\ (y,p)$. Thus  $(V*V')_{\alpha,\beta}=0$ if $(\alpha,\beta)\ne (x,x),\ (p,p)$, $(x,p),\ (p,x)$. Moreover
\begin{alignat*}{2}(V*V')_{x,x}&=V_{x,y}\otimes V'_{y,x}=\rho\otimes\rho*=\rho\otimes\rho,\\
(V*V')_{x,p}&=V_{x,y}\otimes V'_{y,p}=\rho\otimes g\rho*=\rho\otimes g\rho,\end{alignat*}
where $g$ is the non-trivial element in $\mathbb Z/2\mathbb Z$.

Since both $\{(x,x),\ (p,p)\}$ and $\{(x,p),\ (p,x)\}$ are  $F$-orbits in $Y\times Y$, we get
$$V*V'=U_1\oplus U_2,$$
where $U_1$ and $U_2$ are $F$-v.b. on $Y\times Y$ defined by
\begin{alignat*}{2} &U_{1(x,x)}=(V*V')_{x,x}=\rho\otimes\rho,\qquad U_{1(\alpha,\beta)}=0\quad\text{if\ } (\alpha,\beta)\ne (x,x),\ (p,p);\\
&U_{2(x,p)}=(V*V')_{x,p}=\rho\otimes g\rho,\qquad U_{2(\alpha,\beta)}=0\quad\text{if\ } (\alpha,\beta)\ne (x,p),\ (p,x).\end{alignat*}

According to fact 4.2(a), both $\rho\otimes\rho$ and $\rho\otimes g\rho$ can be decomposed into a direct sum of  $\dim\rho$ irreducible representations of $\sotc$. This means that $V*V'=U_1\oplus U_2$ can be decomposed into a direct sum of  $2\dim\rho$ irreducible $F$-vector bundles on $Y\times Y$. By assumption 4.1 and fact 4.2(a), for $\displaystyle t_{w^{-1}}t_w=\sum_u\gamma_{w^{-1},w,u}t_u$, we get that  $\sum_u\gamma_{w^{-1},w,u}=2\dim\rho$. By assumption 4.1 and formula (18) in subsection 4.2, we get $2\dim\rho=2$, hence $\dim\rho=1$. This forces that $\rho$ is the trivial representation of $\sotc.$

Using formula (18) in subsection 4.2 and  assumption 4.1, we must have $\pi(s_1s_2s_1)=U_1$ and $\pi(s_1s_2s_1s_0s_3s_2s_1)=U_2$ since $s_1s_2s_1$ is a distinguished involution.

 Let $v=s_1s_2s_1s_0s_3s_2s_1$. By assumption 4.1, we have $\pi(t_vt_v)=U_2*U_2=U_1$. This contradicts formula (19) in subsection 4.2 which shows that $t_vt_v\ne t_{s_1s_2s_1}$. So this case would not occur.

\medskip

Case 4. The element $x$ is  $F$-stable and $y$ is not $F$-stable. Then  $Fy$ has cardinality 2. Let $Fy=\{y,q\}$. In this case $F_{x,y}=SO_3(\mathbb C)$ and $\rho$ is an irreducible representation of $\sotc$. We have
 $V_{\alpha,\beta}=0$ if $(\alpha,\beta)\ne (x,y),\ (x,q)$, and $V'_{\alpha,\beta}=0$ if $(\alpha,\beta)\ne (y,x),\ (q,x)$. Thus  $(V*V')_{\alpha,\beta}=0$ if $(\alpha,\beta)\ne (x,x)$. Moreover
$$(V*V')_{x,x}=V_{x,y}\otimes V'_{y,x}\oplus V_{x,q}\otimes V'_{q,x}=\rho\otimes\rho^*\oplus g\rho\otimes g\rho^*=\rho\otimes\rho\oplus g\rho\otimes g\rho,$$
where $g$ is the non-trivial element in $\mathbb Z/2\mathbb Z$. Let $\epsilon$ be the sign representation of $\mathbb Z/2\mathbb Z$ and is regarded naturally as representation of $\otc$.

If the highest weight of $\rho$ is positive, then $\rho\otimes\rho\oplus g\rho\otimes g\rho$ can be decomposed into $2\dim\rho$ irreducible representations of $\otc$. This implies that $V*V'$ can be decomposed into a direct sum of  $2\dim\rho$ irreducible $F$-vector bundles on $Y\times Y$. Since $\dim\rho=1+$ highest weight of $\rho$, and the highest weight of $\rho$ is an even number, we get $2\dim\rho\ge 6$. By assumption 4.1 and fact 4.2(a), for $\displaystyle t_{w^{-1}}t_w=\sum_z\gamma_{w^{-1},w,u}t_u$, we get that  $\sum_u\gamma_{w^{-1},w,u}=2\dim\rho\ge 6$. This contradicts  formula (18) in subsection 4.2, which shows $\sum_u\gamma_{w^{-1},w,u}=2.$

If the highest weight of $\rho$ is 0, that is, $\rho$ is the trivial representation of $\sotc$. Then $\rho\otimes\rho\oplus g\rho\otimes g\rho$ is the direct sum of trivial representation of $\otc$ and $\epsilon$ (regarded as representation of $\otc$).  This implies that $V*V'$ can be decomposed into the direct sum of  $2$ irreducible $F$-vector bundles $U$ and $U'$ on $Y\times Y$. We have $U_{x,x}=$ trivial representation of $\otc$ and $U'_{x,x}=\epsilon$, and $U_{\alpha,\beta}=0=U'_{\alpha,\beta}$ if $(\alpha,\beta)\ne (x,x)$.

Using formula (18) in subsection 4.2 and  assumption 4.1, we must have $\pi(s_1s_2s_1)=U$ and $\pi(s_1s_2s_1s_0s_3s_2s_1)=U'$ since $s_1s_2s_1$ is a distinguished involution.

 Let $v=s_1s_2s_1s_0s_3s_2s_1$. By assumption 4.1, we have $\pi(t_vt_v)=U'*U'=U$. This contradicts formula (19) in subsection 4.2 which shows that $t_vt_v\ne t_{s_1s_2s_1}$.

 In conclusion,  this case would not occur.

\medskip

We have shown that all four possibilities would not occur. This implies that one could not find the $F$-set in the assumption 4.1, so that Lusztig's conjecture for the structure of $J_D$ needs modification.

Bezrukavnikov and Ostrik have shown that $J_D$ is isomorphic to some $K_F(Y\times Y)$, where $Y$ is a finite {\it $F$-set of centrally extended points}, see Theorem 4 in [BO]. It would be interesting to describe the isomorphism explicitly, which is useful to compute irreducible representations of affine Hecke algebras of type $\tilde B_3$.

\medskip

In [X1, 5.15] Xi conjectures that there is a ring isomorphism between $J_D$ and $K_F(\mathcal B_u \times \mathcal B_u )$, where $u$ is an element in the unipotent class corresponding to $D$ and $\mathcal B_u$ is the Springer fiber associated with $u$. Since $\mathcal B_u$ is singular in general, it is better to conjecture that there is a ring isomorphism between $J_D$ and $K_F(\mathcal B_u^{\mathbb C^*}\times \mathcal B_u^{\mathbb C^*})$, where the  $\mathbb C^*$-action on $\mathcal B_u$ is defined in [DLP], which commutes with the action of $F$ on $\mathcal B_u$, and the variety of fixed points under the $\mathbb C^*$-action is smooth (see [DLP]).

\bigskip

\noindent{\bf Acknowledgement:} We thank R. Bezrukavnikov for kindly communicating his notice that Lusztig's conjecture describing the summand of J in terms of equivariant
sheaves on the square of a finite set does not hold as stated for the 2-sided cell
corresponding to the unipotent class in Sp6(C) with 3 equal Jordan blocks and for helpful discussions, and  thank G. Lusztig  for helpful comments. Part of the work was done during
YQ's visit to the Academy of Mathematics and Systems Science, Chinese Academy of Sciences. YQ
is very grateful to the AMSS for hospitality and for
financial supports.


\end{document}